\documentclass[reqno, 12pt]{amsart}

\usepackage{amssymb}
\usepackage{latexsym}
\usepackage{epsfig}
\usepackage{epsf}
\usepackage{float}
\usepackage{a4wide}
\usepackage{hyperref}

\newtheorem{theorem}{Theorem}
\newtheorem{lemma}[theorem]{Lemma}
\newtheorem{corollary}[theorem]{Corollary}
\newtheorem{proposition}[theorem]{Proposition}

\def\bR{{\mathbb R}}  
\def\bN{{\mathbb N}}  
\def\bP{{\mathbb P}}    
\def\bZ{{\mathbb Z}}

\def\bE{{\mathbb E}}
\def\bF{{\mathbb F}}

\def\bV{{\mathbb V}}

\def\calA{{\mathcal A}}
\def\calF{{\mathcal F}}

\def\calV{{\mathcal V}}

\def\calG{{\mathcal G}}

\def\calD{{\mathcal D}}
\def\calO{{\mathcal O}}
\def\calP{{\mathcal P}}

\def\A{{\mathbf A}}
\def\B{{\mathbf B}}
\def\C{{\mathbf C}}
\def\D{{\mathbf D}}

\def\I{{\mathbf I}}

\def\0{{\mathbf 0}}
\def\1{{\mathbf 1}}
\def\e{{\mathbf e}}
\def\n{{\mathbf n}}

\def\v{{\mathbf v}}
\def\u{{\mathbf u}}
\def\x{{\mathbf x}}
\def\y{{\mathbf y}}
\def\z{{\mathbf z}}

\def\reff#1{(\ref{#1})}
\def\proofof #1{{\noindent \emph{Proof of #1}.}}
\def\endproof{{\flushright $\square$}} 

\begin{document}

\title[Asymptotics for first-passage times on Delaunay triangulations]{Asymptotics for first-passage times on\\ Delaunay triangulations}
\author{Leandro P. R. Pimentel}
\address{Institute of Mathematics\\
Federal University of Rio de Janeiro\\ 
Brazil.} 
\email{leandro@im.ufrj.br}
\urladdr{}

\keywords{}
\subjclass[2000]{Primary: 60K35; Secondary: 82B}


\begin{abstract}
In this paper we study planar first-passage percolation (FPP) models on random Delaunay triangulations. In \cite{VW90}, Vahidi-Asl and Wierman showed, using subadditivity theory, that the rescaled first-passage time converges to a finite and non negative constant $\mu$. We show a sufficient condition to ensure that $\mu>0$ and derive some upper bounds for fluctuations. Our proofs are based on percolation ideas and on the method of martingales with bounded increments. 
\end{abstract}

\maketitle

\section{Introduction}
\subsection{Brief historical introduction} The classical planar first-passage percolation (FPP) model \cite{HW65} is constructed on the $\bZ^2$ nearest neighbor graph (or lattice) by attaching i.i.d non-negative random variables $\tau_\e$, with common distribution $\bF$, to edges $\e$ on the underlying graph. The \emph{passage time} along a lattice path $\gamma$ is defined by 
\begin{equation}\label{passage}
t(\gamma):=\sum_{\e\in\gamma}\tau_\e\,,
\end{equation}
while the \emph{first passage time} between two vertices $\u$ and $\v$ is defined by 
\begin{equation}\label{firstpassage}
T(\u,\v):=\inf \{t(\gamma)\,:\,\gamma\in\Gamma(\u,\v) \}\,,
\end{equation}
where $\Gamma(\u,\v)$ denotes the set of all lattice paths connecting $\u$ to $\v$. 

In \cite{HW65}, Hammersley and Welsh investigate the asymptotic behavior, as $n\to\infty$, of the first-passage time from $(0,0)$ to $(n,0)$, here denoted by $T_n$. They notice that the sequence $(\bE T_n)_{n\geq 1}$ is subadditive, and hence, $\bE T_n/n$ converges to some $\mu=\mu(\bF)\in[0,\infty]$ (it can be $0$ or $\infty$). They formulate the fundamental question as follows: Under what condition on $\bF$ does $n^{-1}T_n$ converge to $\mu$? It turns out that this question was the motivation for a crucial advance in the theory of subadditive processes that came in 1968 with Kingman's subadditive ergodic theorem \cite{K68}. This theorem allows one to get the almost-sure (and in $L^1$) convergence of $n^{-1}T_n$ to $\mu$ provided $\bE T_1<\infty$, which turns out to be equivalent to $\mu(\bF)<\infty$.

With the law of large numbers for $T_n$ in hands, the subsequent task is to determine what is the right order of $T_n-\bE T_n$: Does the variance of $T_n$ behave like $n^{2\chi}$ for some $\chi\in(0,1)$? Heuristics arguments \cite{KPZ86} indicate  that $\chi=1/3$ for any well behaved passage time distribution $\bF$. However, the only models for which this has been proved are certain last-passage percolation models related to random permutations \cite{BDJ99,J00}. In the lattice first-passage percolation set-up with general i.i.d. passage times, Kesten \cite{K93} developed a remarkable martingale technique which gave the non trivial bound $\chi\leq 1/2$. Here we are mostly concerned with the application of this technique to a random graph version of the FPP model introduced by Vahidi-Asl and Wierman \cite{VW90} \footnote{Howard and Newman  \cite{HN01} also applied this martingale technique in what they called a euclidean first-passage percolation model.}. In our set up, the underlying graph will be random (a Poisson based Delaunay triangulation) and, for each realization of the random graph, the construction of the FPP model will parallel the classical one.

\subsection{The Delaunay FPP model} 
\begin{figure}[htb]
\begin{center}
\includegraphics[width=0.6\textwidth]{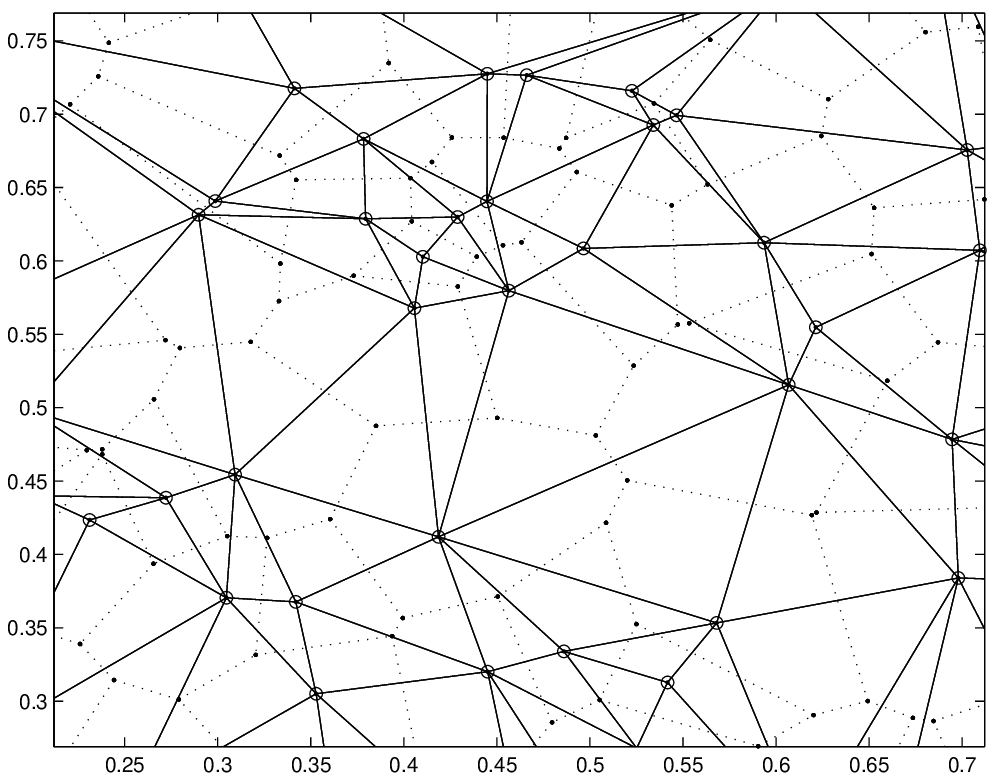}
\end{center}
\caption{The Voronoi Tiling $\calV$ (dotted lines) and the Delaunay Triangulation $\calD$ (solid lines). Most of the work will be on $\calD$.}\label{F:vor2}
\end{figure}

Let $\calP\subseteq\bR^2$ denote the set of points realized in a two-dimensional homogeneous Poisson point process of intensity $1$. To each $\v\in\calP$
corresponds a polygonal region $\C_\v$, named the \emph{Voronoi tile} at $\v$, consisting of points $\x\in\bR^2$ such that $|\x-\v|< |\x-\bar{\v}|$ for all $\bar{\v}\in\calP$, $\bar{\v}\ne\v$. We also denote by $\bar{\C}_\v=\C_\v\cup\partial \C_\v$ the closure of the tile $\C_\v$. The family composed by Voronoi tiles is called the \emph{Voronoi tiling} of the plane based on $\calP$. The \emph{Delaunay triangulation} $\calD$ is the graph where the vertex set is $\calP$ and the edge set  consists of non-oriented pairs $(\v,\bar{\v})$ such that $\C_\v$ and $\C_{\bar\v}$ share a one-dimensional boundary (Figure \ref{F:vor2}). One can see that a.s. each Voronoi tile is a convex and bounded polygon, and the graph $\calD$ is a triangulation of the plane. The \emph{Voronoi tessellation} $\calV$ is the graph where the vertex set consists of vertices of the Voronoi tiles and the edge set is the set of edges of the Voronoi tiles. The edges $\e^*$ of $\calV$ are segments of the perpendicular bisectors of the edges $\e$ of $\calD$. This establishes duality of $\calD$ and $\calV$ as planar graphs.

Each edge $\e\in\calD$ is
independently assigned a nonnegative random variable $\tau_\e$ from a common
distribution $\bF$, which is also independent of the Poisson point process that generates $\calP$. We assume that both $\calP$ and $\tau:=\{\tau_\e\,:\,\e\in\calD\}$ are functions of a configuration $\omega\in\Omega$ and denote by $\bP$ its joint law. The expectation and the variance are denoted by $\bE$ and by $\bV$, respectively. The passage time along a path $\gamma$ on $\calD$ is defined in the same way as in \reff{passage}, while the first-passage time between two vertexes $\v$ and $\v'$ is defined by \reff{firstpassage} where, now, $\Gamma(\v,\v')$ denotes the set of all Delaunay paths connecting $\v$ to $\v'$. A geodesic  connecting $\v$ to $\v'$ is a path which attains the minimum in \reff{firstpassage}: 
\begin{equation}\label{geodesic}
\rho(\v,\v'):=\arg\min\{t(\gamma)\,:\,\gamma\in\Gamma(\v,\v')\}\,. 
\end{equation}
It is known that if $\bF(0)$ is not too close to $1$ then a.s. geodesics do exist for all $\v$ and $\v'$  \cite{P06}. They may not be unique but all of them stay inside an euclidean ball of radius $r=\calO(|\v-\v'|)$. In particular, we almost surely have a finite number of geodesics connecting $\v$ to $\v'$. From now on, $\rho$ will be one of the geodesics that is chosen in some arbitrary way. For each $\x\in\bR^2$ we denote $\v(\x)$ the almost-surely unique point $\v\in\calP$ such that $\x\in\C_{\v}$. For $\x,\y\in\bR^2$ let 
\[
T(\x,\y)\,:=\,T(\v(\x) ,\v(\y)) 
\]
and 
\[
\rho(\x,\y)\,:=\,\rho(\v(\x) ,\v(\y))\,. 
\] 

We recall a fundamental result in the subject, proved by Vahidi-Asl and Wierman \cite{VW90}. As before, for each integer $n\geq 1$, let $T_n:=T(\0,\n)$, where $\0:=(0,0)$ and $\n:=(n,0)$, and 
\begin{equation}\label{timeconstant}  
\mu(\bF)\,:=\,\inf_{n\geq 1}\frac{\bE T_n}{n}\in[0,\infty]\,.
\end{equation}
(Notice that $\mu\leq \bE T_1$.) Assume that $\tau_1,\tau_2,\tau_3$ are independent random variables with common distribution $\bF$: if 
\begin{equation}\label{c1}
\bE(\min_{j=1,2,3}\{\tau_j\})<\infty
\end{equation}
then $\mu(\bF)<\infty$ and, for all unit vectors $\vec{\x}\in\bR^2$ ($|\vec{\x}|=1$),
\begin{equation}\label{e5}
\lim_{n\to\infty}\frac{T(\0,n\vec{\x})}{n}\stackrel{a.s.}{=}\lim_{n\to\infty}\frac{\bE T_n}{n}=\mu(\bF)\,.
\end{equation} 
The main advantage of Poisson-based FPP models is that it provides invariance with respect to all rigid motions, and this  implies that the limit \reff{e5} does not depend on the particular direction $\vec\x$. We note that, under condition \reff{c1}, Vahidi-Asl and Wierman \cite{VW90} proved that $\bE T_1<\infty$, which allows them to use the Kingman's subadditive ergodic theorem. 

\subsection{A percolation threshold and the main result} If $\bF(0)$ is too close to one then there will be, with positive probability, an infinite connected set (or cluster) that contains the origin and is  composed by edges with zero passage times (percolation of zero passage times occurs). In this situation, it is easy to guess that $\mu(\bF)=0$. As we mentioned in the introduction, in this paper we are mostly concerned with Kesten's martingale technique applied to the our random graph FPP model. In order to apply this technique we shall avoid the situation described below.     

The bond percolation model on the Voronoi Tessellation $\calV$, with parameter $p$ and probability law denoted by $\bP^*_p$, is constructed by choosing each edge $\e^*$ in $\calV$ to be open independently with probability $p$. An open path is a path composed by open edges. Let $A_L$ be the event that there exists an open path that crosses the rectangle $[0,3L]\times [0,L]$ in the horizontal direction, i.e.: there exists an open self-avoiding path (no edges intersect) $\gamma^*=(\v_j^*)_{1\leq j\leq k}$ on $\calV$ such that $\v_j^*\in[0,3L]\times [0,L]$ for all $j=2,\dots,k-1$, the line segment connecting $\v_1^*$ to $\v_2^*$ intersects $\{0\}\times[0,L]$ and the line segment connecting $\v_{k-1}^*$ to $\v_k^*$ intersects $\{3L\}\times[0,L]$. Consider the threshold,
\begin{equation}\label{critical}
p_{c}^*:=\inf\left\{p>0\,:\,\liminf_{L\to\infty}\bP_p^*(A_L)=1\right\}\,.
\end{equation}
This threshold should be the same threshold as for bond percolation (existence of an infinite open cluster) in the Voronoi tessellation. It should also be equal to $1-p_c$, where $p_c$ is the threshold for bond percolation in the Delaunay triangulation (duality for planar graphs). However, we will not prove these claims in this paper. For a detailed introduction to percolation models we refer to \cite{RB06,G99}.    

\begin{theorem}\label{t1}
If $\bF(0)<1-p^*_c$ and \reff{c1} holds then
\begin{equation}\label{e1t1}
0<\mu(\bF)<\infty \,. 
\end{equation}
If  \reff{c1} is strengthened to
\begin{equation}\label{c2}
\bE(e^{a\tau})=\int e^{at} d\bF(t)<\infty\mbox{ for some }a>0 \,
\end{equation}
then:
\begin{itemize}
\item Given $\epsilon>0$, for all sufficiently large $n$,
\begin{equation}\label{e2t1}
 \bV(T_n)\leq n^{1+\epsilon}\,;
\end{equation}
\item Given $\kappa\in(1/2,1)$ let $\nu(\kappa)=\frac{4\kappa-2}{5\kappa}$. There exist constants $a_0,a_1>0$ such that for all sufficiently large $n$ and $x\in[1,n^{\kappa}]$
\begin{equation}\label{e3t1}
\bP(|T_n-\bE T_n|\geq xn^{\kappa})\leq a_1e^{-a_0 x^{\nu}}\,.
\end{equation}  
\end{itemize}
\end{theorem}

\subsection{Consequences of Theorem \ref{t1}}
Next we state two corollaries of Theorem \ref{t1}. For the sake of convenience (and space) we only give a sketch of the proofs and leave further details to the reader. The first step is to extend \reff{e3t1} to $T_n-\mu n$.
\begin{corollary}\label{t2}
If $\bF(0)<1-p^*_c$ and \reff{c2} holds, then for all $\kappa\in(1/2,1)$ there exists $a_2>0$ such that for all sufficiently large $n$ and for all $x\in [a_2(\log n)^{1/\nu}, n^{\kappa}]$ ($a_0,a_1,\nu$ are as in Theorem \ref{t1})
\begin{equation}\label{e4t1}
\bP(|T_n-\mu n|\geq xn^{\kappa})\leq a_1e^{-a_0\left(\frac{x}{2}\right)^{\nu}}\,.
\end{equation}
\end{corollary} 
\noindent\emph{Sketch of the proof of Corollary \ref{t2}.} It is clear that Corollary \ref{t2} holds as soon as one proves that, for some constant $a_2>0$,
\begin{equation}\label{e4t2}
\mu n\leq \bE T_n\leq \mu n+\frac{a_2}{2} n^{\kappa} (\log n)^{1/\nu}\,,
\end{equation}
for all sufficiently large $n$. Following \cite{HN01}, the idea to prove \reff{e4t2} is to show that $T_n$ exhibits a weak supper-additivity property:
\begin{equation}\label{supper}
\bE T_{2n}\geq 2\bE T_n - 2c n^{\kappa}(\log n)^{1/\nu}\,,
\end{equation}
where $c>0$ does not depend on $n$. By Lemma 4.2 in \cite{HN01}, \reff{supper} implies \reff{e4t2} . The idea behind the proof of \reff{supper} is as follows. Let $\partial D(n)$ be the boundary of the ball centered at $\0$ and of radius $n$. By translation invariance,
$$\bE T_{2n}\geq 2\bE\min_{\x\in \partial D(n)} T(\0,\x)\,.$$  
On the other hand,   
$$\bE \min_{\x\in \partial D(n)} T(\0,\x)\geq \bE T_n-\bE \max_{\x\in \partial D(n)}\{\bE T(\0,\x)- T(\0,\x)\}\,.$$
The tail bound \reff{e3t1} is strong enough to have that, for all sufficiently large $n$,
$$\bE |\max_{\x\in \partial D(n)}\{\bE T(\0,\x)- T(\0,\x)\}|\leq c n^{\kappa}(\log n)^{1/\nu}\,,$$
(Lemma 4.3 in \cite{HN01}), which yields \reff{supper}. \begin{flushright}\endproof\end{flushright}

\begin{corollary}\label{shape}
Let 
\[
\B_\0 (t):=\{\x\in\bR^2\,:\,\x\in\bar\C_\v\mbox{ where }\v\in\calP\mbox{ and }T(\0,\v)\leq t\}\,
\]  
(recall that $\bar\C_\v$ includes the boundary of the tile $\C_\v$).
Assume that $\bF(0)<1-p^*_c$ and that \reff{c2} holds. Then for all $\kappa\in(1/2,1)$, a.s. there exists $t_0 >0$ such that for all $t>t_0$
\[
(t-t^\kappa)D_\0(1/\mu)\subseteq \B_\0(t)\subseteq (t+t^{\kappa})D_\0(1/\mu)\,
\] 
where $D_\z(r) :=\{ \x\in\bR^2\, : \, |\x-\z|\leq r\}$. 
\end{corollary}
\noindent\emph{Sketch of the proof of Corollary \ref{shape}.} First we claim, without proof, that Corollary \ref{shape} holds as soon as we can prove that, a.s., there exists $M>0$ such that 
\begin{equation}\label{e1shape}
|T(\0,\x)-\mu |\x||\leq |\x|^{\kappa}\,,\mbox{ for all $|\x|\geq M$}\,.
\end{equation}
We sketch the proof of \reff{e1shape} as follows. For each $\x\in\bR^2$ let $\z_\x\in\bZ^2$ be the nearest lattice point to $\x$ ($|\x-\z_\x|\leq 1$). If $|T(\0,\x)-\mu|\x||> |\x|^\kappa$ then 
\begin{equation}\label{summable}
|T(\0,\z_\x)-\mu|\z_\x||> |\z_\x|^\kappa/3\,\,\mbox{  or  }\,\,T_{\z_\x} > |\z_\x|^\kappa/3\,,
\end{equation}
where $T_{\z}:=\max_{\y\in D_{\z}(1)}T(\z,\y)$ (sum and subtract $T(\0,\z_\x)-\mu|\z_x|$ inside the norm of $T(\0,\x)-\mu|\x|$, and then use sub-additivity). Fix $\kappa'\in(1/2,\kappa)$. By Corollary \ref{t2}, the probability that 
$$|T(\0,\z)-\mu|\z||> \frac{|\z|^\kappa}{3}=x|\z|^{\kappa'},\,\,\mbox{ where }\,\,x=\frac{|\z|^{\kappa-\kappa'}}{3}\,,$$
is summable over $\z\in\bZ^2$. On the other hand, the distribution of the random variable $T_{\z}$ does not depend on $\z$ and it has finite moments of all orders\footnote{By Lemma \ref{path}, for each $\y$ with $|\y|\leq 1$, one can always find a path $\gamma_\y$ connecting $\0$ to a point $\y$ with a ``small'' number of steps. Together with \reff{c2}, this yields exponential tail bounds for the random variable $\max_{\y\in D_\0(1)}\sum_{\e\in\gamma_\y}\tau_\e$, that is clearly greater than $T_{\0}$.}. Hence, the probability that $T_{\z}\geq   |\z|^\kappa/3$ is also summable over $\z\in\bZ^2$. By Borel Cantelli's lemma, we then have that \reff{summable} can only happen a.s. for finitely many $\z$'s in $\bZ^2$, which shows \reff{e1shape}.
\begin{flushright}\endproof\end{flushright}

\section{Applying the method of bounded differences via full boxes}\label{pre}
Kesten's idea to control the fluctuations of the first-passage time about its expected value was to represent $T_n-\bE T_n$ as a sum of martingales increments and, after estimating these increments, to apply standard bounds for martingales with bounded increments. First we formulate the abstract set up, following Lemma 5.6 in \cite{HN01}, and then we explain how to use it in our context.
\begin{lemma}\label{lm}  Let $(\Omega,\mathcal{F},\bP)$ be a probability space and let
$\{\mathcal{F}_{k}\}_{k\geq 0}$ be an increasing family
of $\sigma$-algebras of measurable sets. Let $\{M_{k}\}_{k\geq 0}$, $M_{0}=0$,
be a martingale with respect to the filtration $\{\mathcal{F}_{k}\}_{k\geq 0}$
and let $\{U_{k}\}_{k\geq 0}$ be a collection of positive random variables
that are $\mathcal{F}$-measurable. Assume that the increments
$\Delta_{k}=M_{k}-M_{k-1}$ satisfy
\begin{equation}\label{e1lm}
 |\Delta_{k}|\leq z,\mbox{ for some constant }z\geq 0 ,
\end{equation}
 and 
\begin{equation}\label{e2lm}
\bE(\Delta_{k}^{2}\mid \mathcal{F}_{k-1})\leq \bE(U_{k}\mid \mathcal{F}_{k-1}).
\end{equation}
 Assume further that for some finite and positive constants $c'_1,\upsilon,x_{0}$ with $x_{0}\geq z^{2}$, we have that for all $x\geq x_{0}$, 
\begin{equation}\label{e3lm}
 \bP(\sum_{k\geq 1}U_{k}>x)\leq c'_1 e^{-x^{\upsilon}}.
\end{equation}
Then, almost surely, $M=\lim_{k\to\infty}M_{k}$ exists. Moreover, there exist positive and finite constants $c_2,c_3$, whose value does not depend on $z$ and $x_{0}$, such that for all $x\leq x_{0}^{\upsilon}$,
\begin{equation}\label{e4lm}
 \bP(|M|>x\sqrt{x_{0}})\leq c_2e^{-c_3x}.
\end{equation}
\end{lemma}
\proofof{Lemma \ref{lm}} See the proof of Lemma 5.6 in \cite{HN01}.
\begin{flushright}\endproof\end{flushright}
 
In our set up $(\Omega,\mathcal{F},\bP)$ will be the product space over $k\geq 1$ of some probability spaces $(\Omega_k,\mathcal{A}_k,\bP_k)$. Thus, any element of $\Omega$ can be written as $\omega=(\omega_k)_{k\geq 1}$. To any integrable random variable $X$ one can associate the martingale (Doob martingale) 
$$
M_k:=\bE(X\mid \omega_1,\dots,\omega_k)-\bE X
$$ 
with respect to the filtration $\calF_k:=\prod_{1\leq j\leq k}\calA_j$ (we set $\calF_0=\{\emptyset,\Omega\}$ and $M_0 = 0$). To perform some calculations, let us introduce the following notation: if $\omega=(\omega_k)_{k\geq 1}$ and $\sigma=(\sigma_k)_{k\geq 1}\in\Omega$ then we denote
\begin{equation}\label{notation}
[\omega,\sigma]_k :=(\omega_1,\dots,\omega_k,\sigma_{k+1},\sigma_{k+2},\dots)\,.
\end{equation}
Hence,
\[
M_k:=\bE(X\mid \omega_1,\dots,\omega_k)-\bE X=\int X[\omega,\sigma]_k\prod_{j=k+1}^{\infty}d\bP_j(\sigma_j)-\int X(\sigma)d\bP(\sigma)\,.
\]
Since
\[
\int X[\omega,\sigma]_k\prod_{j=k+1}^{\infty}d\bP_j(\sigma_j)=\int X[\omega,\sigma]_k\prod_{j=k}^{\infty}d\bP_j(\sigma_j)\,,
\]
we have 
\[
\Delta_k=\int\Big\{X[\omega,\sigma]_k-X[\omega,\sigma]_{k-1}\Big\}\prod_{j=k}^{\infty}d\bP_j(\sigma_j)\,. 
\]
The simple calculation above shows that, in order to use Lemma \ref{lm} for the Doob martingale, we should control  $X[\omega,\sigma]_k-X[\omega,\sigma]_{k-1}$. This difference is the increment we get when we only change the $k$th coordinate of $[\omega,\sigma]_{k-1}$ from $\sigma_k$ to $\omega_k$.  

Kesten's main point to bound the increment of $X=T_n$ was as follows. For clearness sake, let us assume that our underlying graph is fixed and that each $\omega_k$ represents the passage time attached to an edge $\e_k$. Thus, $T_n[\omega,\sigma]_k-T_n[\omega,\sigma]_{k-1}$ is exactly the increment we get when we change the passage time from $\sigma_k$ to $\omega_k$. If the edge $\e_k$ neither belongs to $\rho_n[\omega,\sigma]_k$ nor to $\rho_n[\omega,\sigma]_{k-1}$ (recall the definition \reff{geodesic} of the geodesic $\rho_n=\rho(\0,\n)$), then  $T_n[\omega,\sigma]_k=T_n[\omega,\sigma]_{k-1}$. Therefore, the edges that have an influence on $T_n-\bE T_n=\sum_{k\geq 1}\Delta_k$ are those caught by $\rho_n$. Furthermore, the martingales increments $\Delta_k$ are uncorrelated and we should expect that the increment by changing the value from $\sigma_k$ to $\omega_k$ is bounded by a function of the pair $(\omega_k,\sigma_k)$. This indicates that the variance of $T_n$ is bounded by a constant times the length of a time minimizing path, which should grow at most linearly in $n$.   

The first problem that appears when applying this method to FPP models is assumption \reff{e1lm}, since by bounding the increment we get a random variable (depending on $(\omega_j,\sigma_j)$) that could attain arbitrarily large values with positive (but small) probability. In the $\bZ^2$ context, this problem was managed by Kesten \cite{K93} by changing the passage times $\tau_\e$ to $\tau_\e(n)=\min\{\tau_e,b\log n\}$, for a suitable constant $b>0$. This truncation allow him to get assumption \reff{e1lm} with $z\sim\log n$. Then he showed that, to get the concentration inequalities for the original model, it suffices to get the same result for the truncated version. In our context we shall deal not only with the random passages times but also with the random graph. This brings new difficulties since a local modification of the point process inside a box $\B\subseteq\bR^2$ can affect the geometry of Voronoi tiles arbitrarily far away from $\B$. To solve this problem, we introduce a modified version of the point process as well as truncated passage times. 
    
\subsection{Full boxes and the modified model}\label{trunc}
For each $\z\in\bZ^2$, $r>0$ and $s\in\{j/2\,:\,j\in\bN\}$ let 
\begin{equation}\label{ebox}
\B_\z^{s,r}:= r\z+[-sr,sr]^2\, .
\end{equation}
The random geometry of the Delaunay Triangulation is controlled through the notion of \emph{full boxes}. Precisely, divide a box $\B$ into thirty-six sub boxes of the same length, say $\B_1,\dots,\B_{36}$. Let $\calP$ be a locally finite subset of the plane. We say that $\B$ is a full box with respect to $\calP$ if all these thirty-six sub boxes have at least one point belonging to $\calP$ . We say that $\Lambda:=(\B_{\z_1}^{1/2,r},\dots,\B_{\z_k}^{1/2,r})$ is a circuit of boxes if $(\z_1,\dots,\z_k)$ is a circuit in $\bZ^2$. Let $\Lambda^{out}$ denote the (topological) interior of the unbounded component of $\bR^2\backslash\cup_{j=1}^k\B_{\z_j}^{1/2,r}$ and $\Lambda^{in}$  denote the union of the interior of the bounded components of $\bR^2\backslash\cup_{j=1}^k\B_{\z_j}^{1/2,r}$. 
\begin{lemma}\label{control}
Assume that  $\Lambda:=(\B_{\z_1}^{1/2,r},\dots,\B_{\z_k}^{1/2,r})$ is a circuit of full boxes with respect to $\calP$. Assume further that $\calP'$ only differs from $\calP$ in $\Lambda^{in}$. Then the Voronoi Tilings $\calV(\calP)$ and $\calV(\calP')$ are the same when restrict to tiles that intersect $\Lambda^{out}$.
\end{lemma}
\proofof{Lemma \ref{control}}  This follows from Lemma 2 of \cite{P06}. The geometrical idea behind it is that the existence of a large Voronoi tile implies the existence of a large empty region, and hence, a Voronoi tile can not cross a region with too many points close by. 
\begin{flushright}\endproof\end{flushright}
 
Now we modify the the Poisson point process such that full boxes occur in an appropriate scale $n^\delta$, for a small $\delta>0$. To be precise, order the points of $\bZ^{2}$ in some arbitrary fashion, say $\bZ^2:=\{\u_1,\u_2, \dots\}$. Let $\delta>0$ be a fixed parameter whose value will depend on $\epsilon>0$ and $\kappa>1/2$ given in Theorem \ref{t1} . Fix $n\geq 1$ and for each $k\geq1$ let 
\[
\B_k^{n}:=\B_{\u_k}^{1/2,n^{\delta}}\,. 
\]
Divide $\B_k^n$ into $36$ sub-boxes (as before) of the same length $n^\delta/6$, say  $\B_{k,1}^n\dots\B_{k,36}^n$. Now we construct the modified point process $\calP(n):=\calP(n,\calP)$ (whose distribution will also depend on $\delta$) by changing the original Poisson point process $\calP$ inside each $\B_{k,j}^n$ (recall that $n$ is fixed), as follows. Let $|\A|$ denote the number of elements belonging to the set $\A$. 
\begin{enumerate}
\item If $1\leq |\B_{k,j}^{n}\cap\calP|\leq n^{2\delta}$ then set $\B_{k,j}^{n}\cap\calP(n):=\B_{k,j}^{n}\cap\calP$;
\item If $|\B_{k,j}^{n}\cap\calP|> n^{2\delta}$ then set $\B_{k,j}^{n}\cap\calP(n)$ by uniformly selecting  $n^{2\delta}$ points from $\B_{k,j}^{n}\cap\calP$.
\item If $|\B_{k,j}^{n}\cap\calP|=0$ then set $\B_{k,j}^{n}\cap\calP(n)$ by adding an extra point  uniformly distributed on $\B_{k,j}^{n}$.
\end{enumerate} 
In a few words, we tile the plane into boxes $\B_k^n$ of size $n^\delta$ and we insist that each tile is a full box, and that no tile contains more than $36 n^{2\delta}$ Poisson points. We make the convention $\calP(\infty)=\calP$ and denote by $\calD(n)$ the Delaunay Triangulation based on $\calP(n)$.

To attach passage times to edges we proceed as follows. Consider a collection  $\tau=\{\tau_{k,j}^{l,m}\,:\,(k,j,l,m)\in \bN^{4}\}$ of i.i.d. random variables. The truncated passage times are given by
\begin{equation}\label{coupled}
\tau_{k,j}^{l,m}(n):=\min\{\tau_{k,j}^{l,m},b\log n\}\,,
\end{equation}
where $b>0$ is a constant whose value will depend on $a>0$ given by \reff{c2}, and we denote $\tau(n)=\{\tau_{k,j}^{l,m}(n)\,:\,(k,j,l,m)\in \bN^{4}\}$ (with the convention $\tau=\tau(\infty)$). We are going to denote by $\v_{k,j}$ the $j$'th vertex in $\mathcal{P}(n) \cap \B_k^n$. For each edge $\e\in\calD(n)$ we must have that $\e=( \v_{k,j},\v_{l,m})$ with either $l> k$ and $j,m\geq 1$, or $l=k$ and $m>j$. (It is just saying either the second end is in a later box, or it is a later point in the same box.) Thus, we set $\tau_\e(n)=\tau_{k,j}^{l,m}(n)$. We would like to stress that the modified model is constructed as function of the original one so that they are naturally coupled.  

From now on, we will work on the product probability space $(\Omega,\mathcal{A},\bP)$ induced by the sequence of independent random elements $\omega_k(n)=\big(\B_k^{n}\cap\calP(n),\tau_k(n)\big),$ for $k\geq 1$, where $\tau_k(n):=\{\tau_{k,j}^{l,m}(n)\,:\,l\geq k,\,\,j,m\geq 1\}$ are the passage times for edges pointing to later points or later boxes. We also use the notation $\omega(n)=\big(\omega_k(n)\big)_{k\geq 1}$ (with $\omega(\infty)=\omega$). The next step is to use Lemma \ref{control} to get an upper bound for the martingale increments in this modified model (recall notation \reff{notation}).
\begin{lemma}\label{l9} Let $I_k$ be the indicator function of the event that there exists $\v\in\rho_n=\rho(\0,\n)$ (recall \reff{geodesic}) such that $\C_\v\cap\B_{\u_k}^{3/2,n^{\delta}}\neq\emptyset$. Then 
\begin{eqnarray}
\nonumber |T_n[\omega(n),\sigma(n)]_k -T_n[\omega(n),\sigma(n)]_{k-1}|&\leq& \\
\nonumber 900\, b \,n^{2\delta}\log n&\times& \max\{I_k[\omega(n),\sigma(n)]_k,I_k[\omega(n),\sigma(n)]_{k-1}\}\,.
\end{eqnarray}
\end{lemma}

\proofof{Lemma \ref{l9}} By the construction of the  modified point process, for each $k\geq 1$, $\Lambda_k:=(\B_j^n\,:\,|\u_j-\u_k|_{\infty}=1)$ is a circuit of full boxes surrounding $\B_k^n$. By Lemma \ref{control}, this prevents any change in $\Lambda_k^{in}=\B_k^n$ affecting the Delaunay Triangulation $\calD(n)$ in the outside region $\Lambda_k^{out}=(\B_{\u_k}^{3/2,n^{\delta}})^c$. In particular, if $I_k[\omega(n),\sigma(n)]_k =0$ then $\rho_n[\omega(n),\sigma(n)]_k$ is also a path in $\calD[\omega(n),\sigma(n)]_{k-1}$ (since $[\omega(n),\sigma(n)]_k$ and $[\omega(n),\sigma(n)]_{k-1}$ only differ inside $\B_k^n$) and thus
\[
T_n[\omega(n),\sigma(n)]_{k-1}\leq t(\rho_n[\omega(n),\sigma(n)]_k)=T_n[\omega(n),\sigma(n)]_k\,.
\]
Analogously, if $I_k[\omega(n),\sigma(n)]_{k-1}=0$ then
\[
T_n[\omega(n),\sigma(n)]_{k}\leq t(\rho_n[\omega(n),\sigma]_{k-1})=T_n[\omega(n),\sigma(n)]_{k-1}\,.
\]
Consequently, if 
$$
\max\{I_k[\omega(n),\sigma(n)]_k,I_k[\omega(n),\sigma(n)]_{k-1}\}=0
$$ 
then 
$$
|T_n[\omega(n),\sigma(n)]_k -T_n[\omega(n),\sigma(n)]_{k-1}|=0\,.
$$

Now assume $I_k[\omega(n),\sigma(n)]_k =1$. Let $\calG_k$ be the subgraph of $\calD(n)$ induced by vertices $\v$ such that $\C_{\v}\cap\B_{\u_k}^{3/2,n^{\delta}}\neq\emptyset$. Then $\calG_k$ is a connected graph. Order the vertices in $\rho_n[\omega(n),\sigma(n)]_k=(\v_j)_{j=1,\dots,m}$,  assume that $\v_1$ and $\v_m$ are not in $\calG_k$, and let $\v_{j_f}$ (resp. $\v_{j_l}$) be the first (resp. the last) vertex of $\rho_n[\omega(n),\sigma(n)]_k$ in $\calG_k[\omega(n),\sigma(n)]_k$. We will bound the increment only in this case. The proof in the other cases follows the same argument with some minors changes. By Lemma \ref{control}, $\rho(\0,\v_f)[\omega(n),\sigma(n)]_k$ and $\rho(\v_l,\n)[\omega(n),\sigma(n)]_k$ are paths in $\calD[\omega(n),\sigma(n)]_{k-1}$. Since $\calG_k[\omega(n),\sigma(n)]_{k-1}$ is connected, one can always connect $\v_f$ to $\v_l$ by a path $\gamma(\v_f,\v_l)$ with no repeated vertices. Let
\[
\psi_n:=\rho(\0,\v_f)\gamma(\v_f,\v_l)\rho(\v_l,\n)\, 
\]
be the concatenation of these three paths. Then, $\psi_n$ is a path in $\calD[\omega(n),\sigma(n)]_{k-1}$ and (recall \reff{passage})
 \begin{eqnarray}
 \nonumber t(\psi_n)[\omega(n),\sigma(n)]_{k-1}&=&T(\0,\v_f)[\omega(n),\sigma(n)]_k \\
  \nonumber&+&t(\gamma(\v_f,\v_l))[\omega(n),\sigma(n)]_{k-1}+T(\v_l,\n)[\omega(n),\sigma(n)]_k\,,
 \end{eqnarray}
which shows that 
\begin{eqnarray}
\nonumber T_n[\omega(n),\sigma(n)]_{k-1}&\leq& t(\psi_n)[\omega(n),\sigma(n)]_{k-1}\\
\nonumber &\leq& T_n[\omega(n),\sigma(n)]_k + t(\gamma(\v_f,\v_l))[\omega(n),\sigma(n)]_{k-1}\,.
\end{eqnarray}
Analogously, if $I_k[\omega(n),\sigma(n)]_{k-1}=1$, then 
\begin{eqnarray}
\nonumber T_n[\omega(n),\sigma(n)]_{k}&\leq& t(\psi_n)[\omega(n),\sigma(n)]_{k}\\
\nonumber &\leq& T_n[\omega(n),\sigma(n)]_{k-1} + t(\gamma(\v_f,\v_l))[\omega(n),\sigma(n)]_{k}\,.
\end{eqnarray}
Therefore, Lemma \ref{l9} holds as soon as we can show that 
\begin{equation}\label{add}
\max\big\{t(\gamma(\v_f,\v_l))[\omega(n),\sigma(n)]_{k}\,\,,t(\gamma(\v_f,\v_l))[\omega(n),\sigma(n)]_{k-1}\big\}\leq 900\,b\,n^{2\delta}\log n\,.
\end{equation}
By the construction of the modified model, each passage time picked up by $\gamma(\v_f,\v_l)$ is bounded by $b\log n$. To bound the size of $\gamma(\v_f,\v_l)$ notice that the graph $\calG_k$ stays within $\B_{\u_k}^{5/2,n^{\delta}}$ (Lemma \ref{control}). In particular, the number of vertices in $\calG_k$ is bounded by the number of points in $\calP(n)\cap\B_{\u_k}^{5/2,n^{\delta}}$, which is at most $25\times36\, n^{2\delta}=900\, n^{2\delta}$. Since $\gamma(\v_f,\v_l)$ has no repeated vertices, its size can not exceed the number of vertices in $\calG_k$, which shows \reff{add}. 
\begin{flushright}\endproof\end{flushright}

\section{Tail bounds for first-passage times and the proof of the main result}
The following propositions control the tail probabilities for first-passage times and time minimizing paths. We would like to stress that they provide upper bounds that hold for all modified models simultaneously. The proof of these propositions is an independent part of the paper and, for clearness sake, they are performed in Section \ref{ap}.

\begin{proposition}\label{tail-1}
If $\bF(0)<1-p^*_c$ then there exist finite and positive constants $c_4,c_5,c_6$ such that for all sufficiently large $n$ (including $n=\infty$) and for all $r\geq 1$,
\begin{equation}\label{coltime}
\bP\big(T_r\big(\omega(n)\big)\leq c_4 r\big)\leq c_5 e^{-c_6 r}\, .
\end{equation}
If \reff{c2} holds, then there exists finite and positive constants $y_1,c_7,c_8$ such that for all  sufficiently large $n$ (including $n=\infty$) and for all $r\geq 1$ and $y\geq y_1$
\begin{equation}\label{coltime1}
\bP\big(T_r\big(\omega(n)\big)\geq yr\big)\leq c_ 7e^{-c_8 yr}\, .
\end{equation}
\end{proposition}

Let $|S|$ denote denote the usual cardinality of a set $S$. For a path $\gamma$ in the Delaunay triangulation we also denote the number of edges by $|\gamma|$. 

\begin{proposition}\label{tail-2}
If $\bF(0)<1-p^*_c$ and $\reff{c2}$ holds then there exist finite and positive constants $y_2,c_9,c_{10},c_{11},c_{12},c_{13}$ such that for all sufficiently large $n$ (including $n=\infty$), for all $r\geq 1$ and for all $y\geq y_2$,
\begin{equation}\label{etail-2-1}
\bP\left(|\rho_r\big(\omega(n)\big)|>yr\right)\leq c_{9}e^{-c_{10}yr}\,,
\end{equation}
and
\begin{equation}\label{etail-2-2}
\bP\left(\rho_r \big(\omega(n)\big)\not\subseteq [-c_{11} r,c_{11} r]^2\right)\leq c_{12}e^{-c_{13} r}\,.
\end{equation}
\end{proposition}

Now we show that the modified model is a good approximation of the original model. For better notation we set
\begin{equation}\label{notmod}
\bar T_n:=T_n\big(\omega(n)\big)\,\,\mbox{ and }\,\,\bar\rho_n:=\rho_n\big(\omega(n)\big)\,.
\end{equation}
We note that the modified model is constructed as a function of the original one, so they are naturally coupled. We also recall that the definition \reff{coupled} of the truncated passage times depends on a positive constant $b>0$. This will allow us to prove the following lemma:  

\begin{lemma}\label{l8}
If $\bF(0)<1-p^*_c$ and \reff{c2} holds then one can choose $b>0$ such that there exist finite constants $b_0,b_1,b_2>0$ such that for all sufficiently large $n$
\begin{equation}\label{e1l8}
\bP\left(|T_n -\bar T_n|> x\right)\leq b_0e^{-b_1 n^{2\delta}}+ 2e^{-\frac{a}{2} x}
\end{equation}
and
\begin{equation}\label{e2l8} 
\bE\left(\{T_n-\bar T_n\}^2\right)\leq b_2\, .
\end{equation}
\end{lemma}  

\proofof{Lemma \ref{l8}} Assume that: (i) $\calP\cap[-2c_{11}n,2c_{11}n]^2=\calP(n)\cap[-2c_{11}n,2c_{11}n]^2$ so that the respective Delaunay triangulation match inside $[-c_{11}n,c_{11}n]^2$ (we are using Lemma \ref{control}). Assume further that: (ii) $\rho_n,\bar\rho_n \subseteq [-c_{11} n,c_{11} n]^2$. Under (i) and (ii), $\rho_n$ is also a possible path for $\bar T_n$ and, since the truncated passage times are smaller than the original ones \reff{coupled}, we must have that $\bar T_n\leq T_n$. On the other hand,  $\bar\rho_n$ is also a possible path for $T_n$, so the sum of the original passage times along this path is at least $T_n$. Hence,
$$0\leq T_n -\bar T_n\leq \sum_{e\in\bar\rho_n}\tau_e\I\{\tau_e > b\log n\}\leq \sum_{\e\in\calD(n)\cap [-c_{11} n,c_{11} n]^2}\tau_\e\I\{\tau_\e > b\log n\} \,.$$
Therefore,
\begin{eqnarray}
\label{e6I8-0}\bP(|\bar T_n-T_n|>x)&\leq& \bP\left(\calP\cap[-2c_{11}n,2c_{11}n]^2\neq\calP(n)\cap[-2c_{11}n,2c_{11}n]^2\right)\\
\label{e6I8}&+&\bP\left(\rho_n \not\subseteq [-c_{11} n,c_{11} n]^2\right)+\bP\left(\bar\rho_n \not\subseteq [-c_{11} n,c_{11} n]^2\right)\\
\label{e6l8-1}&+&\bP\left(\sum_{\e\in\calD(n)\cap[-c_{11} n,c_{11} n]^2}\tau_\e\I\{\tau_\e > b\log n\}> x\right)\,.
\end{eqnarray}
 From standard large deviations results for the Poisson point process, one can see that there exists finite and positive constants $c,c'$ such that right hand side of \reff{e6I8-0} is bounded by $ce^{-c' n^{2\delta}}$. By Proposition \ref{tail-2}, both probabilities in \reff{e6I8} are bounded by $c_{12}e^{-c_{13} n}$, for all sufficiently large $n$. 
 
To deal with \reff{e6l8-1}, we first note that the number of edges in $\calD(n)\cap[-c_{11} n,c_{11} n]^2$ is comparable with the number of points in $\calP(n)\cap[-2c_{11}n,2c_{11}n]^{2}$ (by Euler's relation for planar graphs). By the construction of the modified model, this number is of order $n^2$ (recall that no tile $\B_k^n$ contains more than $36\,n^{2\delta}$ Poissonian points). Thus, $|\calD(n)\cap[-c_{11} n,c_{11} n]^2|\leq c''n^{2}$ for a positive constant $c''$. Now, fix $k\in[1,c''n^{2}]$ and let $\tau_1,\dots,\tau_k$ be i.i.d. random variables with $\bE(e^{a\tau_1})<\infty$ (recall \reff{c2}). By Markov's inequality,
\begin{eqnarray}
\nonumber \bP\left(\sum_{i=1}^{k}\tau_i\I\{\tau_i> b\log n\} > x\right)&\leq& e^{-\frac{a}{2}x}\Big[\bE\left(e^{\frac{a}{2}\tau_1\I\{\frac{a}{2}\tau_1>\frac{ab}{2}\log n\}}\right)\Big]^k\\
\nonumber&=&  e^{-\frac{a}{2}x}\Big[\bE\left(e^{\frac{a}{2}\tau_1\I\{e^{a\tau_1}>n^{\frac{ab}{2}}e^{\frac{a\tau_1}{2}}\}}\right)\Big]^k\,.
\end{eqnarray}
Since 
$$e^{\frac{a}{2}\tau_1\I\{e^{a\tau_1}>n^{\frac{ab}{2}}e^{\frac{a\tau_1}{2}}\}}\leq 1+ e^{\frac{a}{2}\tau_1}\I\{e^{a\tau_1}>n^{\frac{ab}{2}}e^{\frac{a\tau_1}{2}}\}\leq 1 +\frac{e^{a\tau_1}}{n^{\frac{ab}{2}}}\,$$
we get that
$$ \bP\left(\sum_{i=1}^{k}\tau_i\I\{\tau_i> b\log n\} > x\right)\leq e^{-\frac{a}{2}x}\Big[1+\frac{\bE e^{a\tau_1}}{n^{\frac{ab}{2}}}\Big]^k\leq e^{-\frac{a}{2}x}\Big[1+\frac{\bE e^{a\tau_1}}{n^{\frac{ab}{2}}}\Big]^{c''n^2}\leq 2e^{-\frac{a}{2}x} \,,$$
for a sufficiently large $b>0$. Together with the previous bounds, this yields that 
$$\bP(|\bar T_n-T_n|>x)\leq ce^{-c' n^{2\delta}}+2c_{12}e^{-c_{13} n}+2e^{-\frac{a}{2}x}\,,$$
which shows \reff{e1l8} (by choosing suitable $b_0$ and $b_1$). 

By Proposition \ref{tail-1}, we also have that
\begin{equation}\label{e7l8}
\bP(|\bar T_n-T_n|>ny)\leq \bP(\bar T_n>ny/2)+\bP(T_n>ny/2)\leq 2c_7e^{-\frac{c_8}{2}yn}\,,
\end{equation}
if $y>2 y_1$. Together with \reff{e1l8}, this implies that
\begin{eqnarray}
\nonumber\bE\left(\{\bar T_n-T_n\}^2\right)&=&2\int_0^\infty\bP(|\bar T_n-T_n|>z)zdz\\
\nonumber&=&2n^2\int_0^\infty\bP(|\bar T_n-T_n|>ny)ydy\\
\nonumber&=& 2n^2\left[\int_0^{2y_1}\bP(|\bar T_n-T_n|>ny)ydy + \int_{2y_1}^\infty\bP(|\bar T_n-T_n|>ny)ydy\right]\\
\nonumber&\leq&2n^2\left[ \int_0^{2y_1}(b_0e^{-b_1 n^{2\delta}}+ 2e^{-\frac{a}{2} yn})ydy+\int_{2y_1}^\infty2c_7e^{-\frac{c_8}{2}yn}ydy\right]\,,
\end{eqnarray}
which shows \reff{e2l8} 
\begin{flushright}\endproof\end{flushright}

\subsection{Applying Lemma \ref{lm} in the modified model context}
The next step is to apply Lemma \ref{lm} in the modified model context. We consider the martingale 
$$M_k:=\bE\left(\bar T_n\mid\mathcal{F}_k\right)-\bE \bar T_n\,,$$
where the $\sigma$-algebra $\mathcal{F}_k$ is generated by the random elements $\omega_1(n),\dots,\omega_k(n)$ constructed in Section \ref{trunc}. Recall that, in the product space, the increments of this martingale are given by 
\begin{equation}\label{increm}
\Delta_k:=M_k-M_{k-1}=\int\Big\{T_n[\omega(n),\sigma(n)]_k-T_n[\omega(n),\sigma(n)]_{k-1}\Big\}\prod_{j=k}^{\infty}d\bP_j(\sigma_j)\,. 
\end{equation}
By using Lemma \ref{l9}, we will bound this increment with the following variable. Recall the definition of the indicator function $I_k$ (Lemma \ref{l9}) and set $C=(900\,b)^2$.
\begin{lemma}\label{boundincrem-1}
Let 
\begin{equation}\label{boundincrem}
z=\sqrt{C}n^{2\delta}\log n\,\,\mbox{ and }\,\,U_k:= 2C n^{4\delta}(\log n)^2 I_k\big(\omega(n)\big)\,.
\end{equation}
Then
$$|\Delta_k|\leq z\,\,\mbox{ and }\,\,\bE\left(\Delta^2_k\mid\mathcal{F}_{k-1}\right)\leq \bE\left(U_k\mid\mathcal{F}_{k-1}\right)\,.$$
\end{lemma}
\proofof{Lemma \ref{boundincrem-1}} The first inequality follows from Lemma \ref{l9}. By using Lemma \ref{l9} again, and Schwarz' inequality, one has that
\begin{eqnarray}
\nonumber&&\bE(\Delta_k^2\mid\calF_{k-1})=\int\left\{\int\{T_n[\omega(n),\sigma(n)]_k-T_n[\omega(n),\sigma(n)]_{k-1}\}\prod_{j=k}^\infty d\bP_j(\sigma_j)\right\}^2 d\bP_k(\omega_k)\\
\nonumber&&\leq C n^{4\delta}(\log n)^2\int\left\{\int\max\{I_k[\omega(n),\sigma(n)]_k,I_k[\omega(n),\sigma(n)]_{k-1}\}\prod_{j=k}^\infty d\bP_j(\sigma_j)\right\}^2 d\bP_k(\omega_k)\\
\label{eboundincrem-1}&&\leq C n^{4\delta}(\log n)^2\int\int\max\{I_k[\omega(n),\sigma(n)]_k,I_k[\omega(n),\sigma(n)]_{k-1}\}\prod_{j=k}^\infty d\bP_j(\sigma_j) d\bP_k(\omega_k)\,
\end{eqnarray}
where in the last step we also used that the square of the max of indicator functions is just the max of indicator functions. Notice also that
\begin{equation}\label{eboundincrem-2}
\max\{I_k[\omega(n),\sigma(n)]_k,I_k[\omega(n),\sigma(n)]_{k-1}\}\leq I_k[\omega(n),\sigma(n)]_k+I_k[\omega(n),\sigma(n)]_{k-1}\,,
\end{equation}
and that
\begin{eqnarray}
\nonumber\bE\left(I_k(\omega(n))\mid\mathcal{F}_{k-1}\right)&=&\int\int I_k[\omega(n),\sigma(n)]_{k-1}\prod_{j=k}^\infty d\bP_j(\sigma_j)\bP_k(\omega_k)\\
\nonumber&=&\int\int I_k[\omega(n),\sigma(n)]_{k}\prod_{j=k}^\infty d\bP_j(\sigma_j) d\bP_k(\omega_k)\,.
\end{eqnarray}
Together with \reff{eboundincrem-1} and \reff{eboundincrem-2}, this proves the second inequality. 
\begin{flushright}\endproof\end{flushright}
 Now we will see that $U_k$ satisfies \reff{e3lm} (with $x_0=n^{1+5\delta}$).
\begin{lemma}\label{l10}
If $\bF(0)<1-p^*_c$ and \reff{c2} holds then for all sufficiently large $n$
\begin{equation}\label{e1l10}
\bE(\sum_{k=1}^{\infty}U_k)\leq n^{1+5\delta}\,.
\end{equation}
Further, if $\delta\in(0,1/3)$, for all sufficiently large $n$ and $x\geq n^{1+5\delta}$,
\begin{equation}\label{e2l10}
\bP(\sum_{k=1}^{\infty}U_k > x)\leq c_9 e^{-x^{1/2}}\, ,
\end{equation}
($c_9$ as in Proposition \ref{tail-2}).
\end{lemma}

\proofof{Lemma \ref{l10}} If $\v\in\B_{\u_l}^{1/2,n^\delta}$ then $\C_\v\subseteq\B_{\u_l}^{3/2,n^\delta}$ (Lemma \ref{control}), and there are $25$ boxes $\B_{\u}^{3/2,n^\delta}$ that intersect this set.  Hence, $\sum_{k=1}^{\infty}I_k$ is at most $25|\rho_n\big(\omega(n)\big)|$ and therefore,
\begin{equation}\label{e1pl10}
\sum_{k=1}^{\infty}U_k =2C n^{4\delta}(\log n)^2\sum_{k=1}^{\infty}I_k\leq 50Cn^{4\delta}(\log n)^2|\rho_n\big(\omega(n)\big)|\,.
\end{equation}
On the other hand, by Proposition \ref{tail-2}, \reff{etail-2-1}, there exist constants $c_9,c_{10}$ and $y_2$ such that  
$$\bE|\rho_n\big(\omega(n)\big)|\leq \left(y_2+\frac{c_9}{c_{10}}\right)n\,.$$
Together with \reff{e1pl10}, this shows \reff{e1l10} for sufficiently large $n$. Now, if $x\geq n^{1+5\delta}$ then 
$$\frac{x}{50Cn^{4\delta}(\log n)^2}\geq y_2 n\,\,\mbox{ and }\,\,c_{10}\frac{x}{50Cn^{4\delta}(\log n)^2}\geq x^{1/2}\,,$$
for large enough $n$, provided $\delta\in(0,1/3)$. Thus, again by Proposition \ref{tail-2}, \reff{etail-2-1},
\begin{eqnarray}
\nonumber\bP\left(\sum_{k=1}^{\infty}U_k>x\right)&\leq& \bP\left(|\rho_n\big(\omega(n)\big)|>\frac{x}{50C n^{4\delta}(\log n)^2}\right)\\
\nonumber&\leq &c_{9}e^{-c_{10}\frac{x}{50Cn^{4\delta}(\log n)^2}}\leq c_9 e^{-x^{1/2}}\,.
\end{eqnarray}
which implies \reff{e2l10}.  
\begin{flushright}\endproof\end{flushright}

\subsection{Proof of Theorem \ref{t1}} 
We note that \reff{e1t1} follows from Theorem 1 and Corollary 1 in \cite{P06}. By Lemma \ref{boundincrem-1} and Lemma \ref{l10}, with $M_k=\bE(\bar T_n\mid\calF_k)$,
\begin{equation}\label{e2pt1}
\bV(\bar T_n) =\sum_{k=1}^{\infty}\bE(\Delta_k^2)\leq \bE(\sum_{k=1}^{\infty}U_k)\leq n^{1+5 \delta}\,.
\end{equation}
On the other hand, by Lemma \ref{l8},
$$\bV \left(T_n -\bar T_n\right)=\bE\left(\{T_n-\bar T_n\}^2\right)-\left(\bE\{T_n-\bar T_n\}\right)^2\leq \bE\left(\{T_n-\bar T_n\}^2\right)\leq b_2\,.$$
By \reff{e2pt1}, this proves \reff{e2t1}. 

Now, by combining Lemma \ref{boundincrem-1} together with Lemma \ref{l10}, one verifies assumptions  \reff{e1lm}, \reff{e2lm} and \reff{e3lm} of Lemma \ref{lm}, with $x_0=n^{1+5\delta}$ and $\upsilon=1/2$. Therefore, for all sufficiently large $n$ and $x\leq n^{\frac{1+5\delta}{2}}$,
\begin{equation}\label{e0pt1}
\bP(|\bar{T}_n-\bE\bar{T}_n|>xn^{\frac{1+5\delta}{2}})\leq c_2 e^{-c_3 x}\,.
\end{equation} 

To finish the proof of the theorem, for each $\kappa\in(1/2,4/3)$ let $\delta=(2\kappa-1)/5\in (0,1/3)$. By \reff{e0pt1}, for sufficiently large $n$ and $x\in[1, n^{\kappa}]$ (sum and subtract $\bar{T}_n-\bE\bar{T}_n$)
\begin{eqnarray}
\nonumber\bP(|T_n-\bE T_n|>xn^{\kappa})&\leq&\bP\Big(\,|T_n-\bar{T}_n|>\frac{xn^{\kappa}}{3}\,\Big)\\
\nonumber&+&\bP\Big(\,|\bar{T}_n-\bE\bar{T}_n|>\frac{xn^{\kappa}}{3}\,\Big)\\
\label{e1pt1} &\leq&\bP\Big(\,|T_n-\bar{T}_n|>\frac{xn^{\kappa}}{3}\,\Big)+c_2e^{-c_3 x/3}\,,
\end{eqnarray}
where, for the first step, we also used that 
$$|\bE T_n-\bE\bar T_n|\leq\bE|T_n-\bar T_n|\leq(\bE(T_n-\bar T_n)^2)^{1/2}\leq\sqrt{b_2}\leq \frac{xn^{\kappa}}{3}\,.$$
By Lemma \ref{l8} 
\begin{equation}\label{e3pt1}
\bP\Big(\,|T_n-\bar{T}_n|>\frac{xn^{\kappa}}{3}\,\Big)\leq b_0e^{-b_1n^{2\delta}}+2e^{-a\frac{xn^{\kappa}}{6}}=b_0e^{-b_1n^{(4\kappa-2)/5}}+2e^{-a\frac{xn^{\kappa}}{6}}\,.
\end{equation}
Let $\nu(\kappa)=\frac{4\kappa-2}{5\kappa}>0$. Thus, if $x\leq n^\kappa$ then $x^\nu\leq n^{(4\kappa-2)/5}$ and hence
$$b_0e^{-b_1n^{(4\kappa-2)/5}}+2e^{-a\frac{xn^{\kappa}}{6}}\leq b_0e^{-b_1x^{\nu}}+2e^{-a\frac{xn^{\kappa}}{6}}\,.$$
Together with \reff{e1pt1} and \reff{e3pt1}, this proves \reff{e3t1} (for suitable constants $a_0$ and $a_1$).

\section{Appendix}\label{ap}
We start by stating some results from \cite{P10} that will be the starting point for proving Proposition \ref{tail-1} and Proposition \ref{tail-2}. These results concern some geometrical aspects of self-avoiding paths on Delaunay triangulations as follows: Let $\Gamma_{\geq r}(n)$ (resp., $\Gamma_{\leq r}(n)$) be the set of all self-avoiding paths $\gamma$ in $\calD(n)$, starting at $\v(\0)\in\calP(n)$ (the nearest point to $\0$) and of size $|\gamma|\geq r$ (resp.,  $|\gamma|\leq r$). For each $C\subseteq\bR^2$ let
$$\A(C):=\left\{\z\in\bZ^d\,\,:\,\,\B_\z^{1/2,L}\cap C\neq\emptyset\,\right\}\,.$$
If $\gamma$ is a self-avoiding path in $\calD(n)$, let 
$$\A(\gamma):=\A\left(\cup_{\v\in\gamma}C_\v\right)\,.$$
\begin{lemma}\label{self}
There exist constants $b_3,b_4,b_5\in(0,\infty)$, that do not depend on $n\geq 1$, such that for all $r\geq1$ 
\begin{equation}\label{e1l1}
\mbox{ if $\,r\geq b_3 s$ then }\,\,\bP\left( \min_{\gamma\in\Gamma_{\geq r}(n)}|\A(\gamma)|\leq s\right)< e^{-\frac{r}{2}}\,;
\end{equation}
and
\begin{equation}\label{e2c1}
\mbox{ if $\,s\geq b_4 r\,$  then }\,\,\bP\left(\max_{\gamma\in\Gamma_{\leq r}(n)}|\A(\gamma)|\geq s\right)< e^{-b_5 s}\,;
\end{equation}
 \end{lemma}
 \proofof{Lemma \ref{self}} See Corollary 2 and Theorem 12 in \cite{P10}.
  \begin{flushright}\endproof\end{flushright}

\begin{lemma}\label{path}
For each $\x\in\bR^2$ there exists almost surely a Delaunay self avoiding path $\gamma(\0,\x)$ that connects $\v_\0$ to $\v_\x$ and only uses vertices of Voronoi tiles that intersect the line segment $[\0,\x]$. Further, there exist constants $b_6,b_7\in(0,\infty)$, that do not depend on $n\geq 1$, such that for all $r\geq1$,
\begin{equation}\label{e1c2}
\mbox{ if $\,r\geq b_6 s\,$  then }\,\,\bP\left(\max_{\x:\|\x\|_2\leq s}|\gamma(\0,\x)|\geq r\right)< 2e^{-b_7 r}\, .
\end{equation}
\end{lemma} 
\proofof{Lemma \ref{path}} See Corollary 4 and Theorem 12 in \cite{P10}.
\begin{flushright}\endproof\end{flushright}

\begin{lemma}\label{reward}
Assume that $\{\tau_\e\,:\,\e\in\D\}$ is a collection of i.i.d. Bernoulli random variables with $\bP(\tau_\e=0)<1-p_c^*$. Then there exist finite constants $n_0,b_8,b_{9}>0$, that do not depend on $n\geq 1$, such that for for all $n\geq n_0$ 
$$\mbox{ if $\,r\geq b_8 s\,$  then }\bP\left(\min_{\gamma\in\Gamma_{\geq r}(n)}\sum_{\e\in\gamma}\tau_\e\leq s\right)\leq 3 e^{-b_{9} s}\,.$$
\end{lemma}

\proofof{Lemma \ref{reward}} See Theorem 13 in \cite{P10}.
 \begin{flushright}\endproof\end{flushright}

\subsection{Proof of Proposition \ref{tail-1}} Consider the geodesic $\rho_r=\rho(\0,(r,0))$ connecting $\0$ to $(r,0)$. Then $\A(\rho_r)\geq r$ and 
\begin{eqnarray}
\nonumber\bP\left(T_r(\omega(n))<s\right)&\leq& \bP\left(|\rho_r(\omega(n))|< t\right)+\bP\left(|\rho_r(\omega(n))|\geq t\,\mbox{ and } T_r(\omega(n))<s\right)\\
\label{pfprop1} &\leq& \bP\left(\max_{\gamma\in\Gamma_{\leq t}(n)}\#\A(\gamma)\geq r\right)+\bP\left(\min_{\gamma\in\Gamma_{\geq t}(n)}\sum_{\e\in\gamma}\tau_\e(n)<s\right)\,.
\end{eqnarray}
 for any $r,s,t\geq 0$. 
 
By Lemma \ref{self}, if $t=r/b_4$ then 
 $$\bP\left(\max_{\gamma\in\Gamma_{\leq t}(n)}\#\A(\gamma)\geq r\right)=\bP\left(\max_{\gamma\in\Gamma_{\leq r/b_4}(n)}\#\A(\gamma)\geq r\right)\leq e^{-b_5 r}\,.$$
To handle with the other term, fix $\epsilon>0$, and define the auxiliary random variables $\tau_\e^{\epsilon}:=1\{\tau_\e> \epsilon\}$. Then $\{\tau_\e^{\epsilon}\,:\,\e\in\calD\}$ is a collection of i.i.d. Bernoulli random variables such that $\epsilon\tau_\e^\epsilon \leq \tau_\e(n)$ and $\bP\left( \tau_\e^\epsilon=0\right)=\bF(\epsilon)$. Since we have assumed that $\bF(0)<1-p_c^*$, we can choose $\epsilon_0>0$ such that $\bF(\epsilon_0)<1-p_c^*$. Notice that this $\epsilon_0$ can be chosen sufficiently small in order that the distribution of $\tau_\e^\epsilon$ does not dependent on $n\geq 1$. Then
 $$\bP\left(\min_{\gamma\in\Gamma_{\geq t}(n)}\sum_{\e\in\gamma}\tau_\e(n)<s\right)\leq \bP\left(\min_{\gamma\in\Gamma_{\geq t}(n)}\sum_{\e\in\gamma}\tau_\e^{\epsilon_0}<\frac{s}{\epsilon_0}\right)\,.$$
By Lemma \ref{reward}, if $s=\epsilon_0r/(b_4b_8)=\epsilon_0 t/b_8$ (recall that $t=r/b_4$) then 
$$\bP\left(\min_{\gamma\in\Gamma_{\geq t}(n)}\sum_{\e\in\gamma}\tau_\e^{\epsilon_0}<\frac{s}{\epsilon_0}\right)=\bP\left(\min_{\gamma\in\Gamma_{\geq r/b_4}(n)}\sum_{\e\in\gamma}\tau_\e^{\epsilon_0}<\frac{r}{b_4b_8}\right)\leq 3e^{-\frac{b_{9}}{b_4b_8} r}\,,$$
and hence
$$\bP\left(\min_{\gamma\in\Gamma_{\geq r/b_4}(n)}\sum_{\e\in\gamma}\tau_\e(n)<\epsilon_0\frac{r}{b_4b_8}\right)\leq  3e^{-\frac{b_{9}}{b_4b_8} r}\,.$$
Together with \reff{pfprop1}, the above inequalities yield 
$$\bP\left(T_r(\omega(n))<\epsilon_0\frac{r}{b_4b_8}\right)\leq e^{-b_5 r}+3e^{-\frac{b_{9}}{b_4b_8} r}\,,$$
which implies \reff{coltime}. 

To prove \reff{coltime1} we use the path $\gamma_r=\gamma(\0,(r,0))$ given by Lemma \ref{path}. Thus, if $t\geq b_6r$ then
 $$ \bP\left(|\gamma_r|>t\right)\leq 2e^{-b_7t}\,.$$
Notice that  
$$T_r(\omega(n))\leq \sum_{\e\in\gamma_r}\tau_\e(n)\leq\sum_{\e\in\gamma_r} \tau_\e\,,$$
and hence,
\begin{equation}\label{pf1prop1}
\bP\left(T_r(\omega(n))>s\right)\leq \bP\left(|\gamma_r|>t\right)+\bP\left(|\gamma_r|\leq t\mbox{ and }\sum_{\e\in\gamma_r}\tau_\e>s\right)\,,
 \end{equation} 
 for any $r,s,t\geq 0$. The passage times $\tau_\e$ are independent of $\gamma_r$, which is only a function of $\calP(n)$. Therefore, under \reff{c2}, Markov's inequality implies that 
\begin{eqnarray}
\nonumber\bP\left(|\gamma_r|\leq t\mbox{ and }\sum_{\e\in\gamma_r}\tau_\e>s\right)&=&\sum_{k=1}^t\bP\left(\sum_{j=1}^k\tau_j>s\right)\bP\left(|\gamma_r|=k\right)\\
\nonumber&\leq&\sum_{k=1}^t e^{-as}\left\{\bE e^{a\tau_1}\right\}^k\bP\left(|\gamma_r|=k\right)\\
\nonumber&\leq&e^{-as}\left\{\bE e^{a\tau_1}\right\}^t\leq e^{-\frac{a}{2}s}\,,
\end{eqnarray}
where $t=\frac{a}{2\log\bE e^{a\tau_1}}s$ ($\tau_j$ for $j\geq 1$ represent i.i.d. copies of $\tau_\e$). Together with \reff{pf1prop1} (and the bound on the size of $\gamma_r$), this shows that
 $$\bP\left(T_r(\omega(n))>s\right)\leq 2e^{-b_7t}+ e^{-\frac{a}{2}s} = 2e^{-\frac{b_7a}{2\log\bE e^{a\tau_1}}s}+ e^{-\frac{a}{2}s}\,,$$
as soon as $s\geq \frac{2\log\bE e^{a\tau_1}b_6}{a}r$, which proves \reff{coltime1} (take $s=yr$).

\subsection{Proof of Proposition \ref{tail-2}}
Notice that 
\begin{eqnarray}
\nonumber\bP\left(|\rho_r(\omega(n))|>s\right)&\leq& \bP\left(T_r(\omega(n))>t\right)+ \bP\left(T_r(\omega(n))\leq t\mbox{ and }|\rho_r(\omega(n))|>s\right)\\
\label{pf-etail-2-1}&\leq&\bP\left(T_r(\omega(n))>t\right)+ \bP\left(\min_{\gamma\in\Gamma_{\geq s}(n)}\sum_{\e\in\gamma}\tau_\e(n)\leq t\right)\,,
\end{eqnarray}
for any $r,s,t\geq 0$. By repeating the argument used in the proof of \reff{coltime}, if $t=\epsilon_0 s/b_8$ then
$$ \bP\left(\min_{\gamma\in\Gamma_{\geq s}(n)}\sum_{\e\in\gamma}\tau_\e(n)\leq t\right)\leq 3e^{-\frac{b_9}{b_8}s}\,.$$
On the other hand, by \reff{coltime1}, if $t\geq  y_1 r$ then
$$\bP\left(T_r(\omega(n))>t\right)\leq c_7 e^{-c_8 t}=c_7 e^{-c_8 \frac{\epsilon_0}{b_8}s} \,.$$
Combining the last two inequalities together with \reff{pf-etail-2-1}, one has that
$$\bP\left(|\rho_r(\omega(n))|>s\right)\leq c_7 e^{-c_8 \frac{\epsilon_0}{b_8}s} + 3e^{-\frac{b_9}{b_8}s}\,,$$
as soon as $s\geq \frac{y_1b_8}{\epsilon_0}r$, which shows \reff{etail-2-1}.

Now, if $\rho_r\not\subseteq[-s,s]^2$ then $|\A(\rho_r)|\geq s$. Together with \reff{etail-2-1} and Lemma \ref{self}, this implies that
\begin{eqnarray}
\nonumber\bP\left(\rho_r(\omega(n))\not\subseteq[-b_4y_2r,b_4y_2r]^2\right)&\leq &\bP\left(|\rho_r(\omega(n))|>y_2 r\right)\\
\nonumber&+&\bP\left(|\rho_r(\omega(n))|\leq y_2 r\mbox{ and }\rho_r(\omega(n))\not\subseteq[-b_4y_1r,b_4y_2r]^2\right)\\
\nonumber&\leq&\bP\left(|\rho_r(\omega(n))|>y_2 r\right)+\bP\left(\max_{\gamma\in\Gamma_{\leq y_2 r}(n)}\#\A(\gamma)\geq b_4 y_2 r\right)\\
\nonumber&\leq& c_9 e^{-c_{10}y_2 r}+e^{-b_5b_4y_2r}\,,
\end{eqnarray}
and the proof of \reff{etail-2-2} is complete.
\newline

\paragraph*{\bf Acknowledgments} 
Part of this work was develop during my doctoral studies at Impa and I would like to
thank my adviser, Prof. Vladas Sidoravicius, for his dedication and encouragement during this period. I also thank the whole administrative staff of IMPA for their assistance and CNPQ for financing my doctoral studies, without which this work would have not been possible. This work also amends some mistaken passages of my thesis \cite{P04} and, for this reason, I should also express gratitude for the assistance given by Prof. Thomas Mountford. Furthermore, we are grateful for careful reading and the many helpful suggestions by anonymous referees.

\end{document}